\documentclass[11pt,reqno]{amsart}
\usepackage{enumerate}
\usepackage[margin=1in]{geometry}
\usepackage{amsfonts}
\usepackage{amssymb}
\usepackage{amsmath}
\usepackage{amsthm}
\usepackage[cmtip,all]{xy}
\usepackage{microtype}
\usepackage{mathtools}
\usepackage{cite}
\usepackage[hidelinks]{hyperref}

\newcommand\divides{\;\text{ divides }\;}
\def\p#1{\left( #1 \right)}
\def\set#1{\left\{ #1 \right\}}
\def\cyc#1{\langle #1 \rangle}
\def\abs#1{\left| #1 \right|}
\def\Gal{\operatorname{Gal}}

\def\F{\mathbb{F}}
\def\N{\mathbb{N}}
\def\Z{\mathbb{Z}}
\def\Q{\mathbb{Q}}

\def\End{\operatorname{End}}
\def\Mat{\operatorname{Mat}}
\def\GL{\operatorname{GL}}
\def\GSp{\operatorname{GSp}}
\def\Sp{\operatorname{Sp}}
\def\PGSp{\operatorname{PGSp}}
\def\PSp{\operatorname{PSp}}
\def\mult{\operatorname{mult}}

\def\normal{\trianglelefteq}
\def\rad{\operatorname{rad}}
\def\tors{\text{tors}}
\def\lcm{\operatorname{lcm}}

\newcommand{\xrightarrowdbl}[2][]{
  \xrightarrow[#1]{#2}\mathrel{\mkern-14mu}\rightarrow
}

\newtheorem{theorem}{Theorem}[section]
\newtheorem{lemma}[theorem]{Lemma}
\newtheorem{proposition}[theorem]{Proposition}
\newtheorem{corollary}[theorem]{Corollary}

\theoremstyle{definition}

\theoremstyle{remark}
\newtheorem{remark}[theorem]{Remark}

\numberwithin{equation}{section}

\begin{document}

\title[A bound for the image conductor]{A bound for the image conductor of a principally polarized abelian variety with open Galois image}

\author{Jacob Mayle}

\address{Department of Mathematics, Statistics, and Computer Science, University of Illinois Chicago, 851 S. Morgan Street, Chicago, IL 60607}

\email{jmayle2@uic.edu}


\begin{abstract} Let $A$ be a principally polarized abelian variety of dimension $g$ over a number field $K$. Assume that the image of the adelic Galois representation of $A$ is an open subgroup of $\GSp_{2g}(\hat{\Z})$. Then there exists a positive integer $m$ so that the Galois image of $A$ is the full preimage of its reduction modulo $m$. The least $m$ with this property, denoted $m_A$, is called the \textit{image conductor} (also called the \textit{level}) of $A$. Jones recently established an upper bound for $m_A$, in terms of standard invariants of $A$, in the case that $A$ is an elliptic curve without complex multiplication. In this paper, we generalize the aforementioned result to provide an analogous bound in arbitrary dimension. 
\end{abstract}

\maketitle

\section{Introduction}
Let $A$ be a principally polarized abelian variety of dimension $g$ over a number field $K$. Let $T(A) \coloneqq \varprojlim A[m]$ denote the adelic Tate module of $A$. The \textit{adelic Galois representation} of $A$ is a continuous homomorphism of profinite groups 
\[ \rho_{A}: G_K \to \GSp_{2g}(\hat{\Z}) \] 
that encodes the action of $G_K \coloneqq \Gal(\overline{K}/K)$ on $T(A)$.

The image of $\rho_A$ is called the \textit{Galois image} of $A$ and, in many cases, is known to be an open subgroup of $\GSp_{2g}(\hat{\Z})$. For instance, Serre established that this is so for elliptic curves without complex multiplication in his celebrated 1972 open image theorem \cite{MR387283}. Serre later generalized his result to certain higher dimensions. 

\begin{theorem}[Serre, 1986 \cite{MR3185222}] Let $A$ be a principally polarized abelian variety of dimension $g$ over a number field $K$. If $\End_{\overline{K}}(A) = \Z$ and $g = 2,6$, or is odd, then  $\rho_{A}(G_K) \subseteq \GSp_{2g}(\hat{\Z})$ is an open subgroup.
\end{theorem}

Due to an example of Mumford \cite[\S 4]{MR248146}, it is known that the above result does not generalize to arbitrary dimension without further hypotheses. In 2011 \cite{MR2820155}, Hall gave a \textit{sufficient} condition for a principally polarized abelian variety of arbitrary dimension to have open Galois image. Kowalski proved, as a consequence, that almost all Jacobians of hyperelliptic curves (in a suitable sense) have open Galois image  \cite[Appendix]{MR2820155}.

Assume that $A$ has open Galois image. For each positive integer $m$, we let 
\[ \bar{\pi}_m: \GSp_{2g}(\hat{\Z}) \twoheadrightarrow \GSp_{2g}(\Z/m\Z) \] 
be the natural projection map. The collection $\{\ker \bar{\pi}_m\}_{m=1}^\infty$ is a neighborhood basis for the identity of $\GSp_{2g}(\hat{\Z})$. Since $\rho_A(G_K) \subseteq \GSp_{2g}(\hat{\Z})$ is an open subgroup, there exists an $m$ so that $\ker \bar{\pi}_m \subseteq \rho_A(G_K)$. The least $m$ with this property is the \textit{image conductor} of $A$, and is denoted by $m_A$. An important observation is that the Galois image of $A$ is the full preimage of the finite group $\bar{\pi}_{m_A}(\rho_A(G_K))$, as we shall discuss in \S\ref{split-stable}.

In a recent paper \cite{MR4190460}, Jones established an upper bound for $m_A$, in terms of standard invariants of $A$, in the case that $A$ is an elliptic curve without complex multiplication. Further, he remarked that his techniques should be able to be extended to prove an analogous result for principally polarized abelian varieties of arbitrary dimension. In this paper, we do precisely that, proving the following.

\begin{theorem} \label{main-thrm} Let $A$ be a principally polarized abelian  variety of dimension $g$ over a number field $K$ and assume that the image of the adelic Galois representation $\rho_{A}: G_K \to \GSp_{2g}(\hat{\Z})$ is open in $\GSp_{2g}(\hat{\Z})$. Then
\[ m_{A} \leq 2 \cdot \mathcal{B}_A \cdot \left[ \GSp_{2g}(\hat{\Z}) : \rho_{A}(G_K) \right]\!, \]
where $m_A$ denotes the image conductor of $A$ and $\mathcal{B}_{A}$ is the square-free constant, depending on $A$, that is defined to be the product of each prime number $\ell \in \N$ that satisfies  at least one of the following conditions:

\begin{enumerate}
\item $K/\Q$ is ramified at $\ell$; 
\item $A$ has bad reduction at some prime ideal of $\mathcal{O}_K$ that lies over $\ell$; or
\item $\ell = 2$, in the case that $g = 2$.
\end{enumerate}
\end{theorem}

\begin{remark} We now consider sharpness of the bound in Theorem \ref{main-thrm} when $g = 2$. Let $A$ be the Jacobian of a genus 2 curve $C/\Q$. Let $\Delta$ denote the discriminant of $C$. Write $\Delta_{\text{sf}}$ to denote the square-free part of $\Delta$. It follows similarly as in the case of elliptic curves \cite[\S2.4]{MR4041152} that
\begin{equation} \rho_A(G_\Q) \subseteq \set{\gamma \in \GSp_{4}(\hat{\Z}) : \epsilon(\gamma) = \chi_A(\gamma)} \label{epsilon-chi} \end{equation}
where $\epsilon$ and $\chi_A$ are defined as follows: The character $\epsilon$ is the map
\[ \epsilon:\GSp_4(\hat{\Z}) \to \GSp_4(\Z/2\Z) \overset{\sim}{\to} S_6 \to \set{\pm 1} \]
given by projection modulo $2$, followed by the signature character on the symmetric group $S_6$. For the character $\chi_A$, first define the constant
\[ d_A = \begin{cases}
\Delta_{\text{sf}} & \Delta_{\text{sf}} \equiv 1 \pmod{4}, \\
4\Delta_{\text{sf}} & \text{otherwise.}
\end{cases} \]
Now $\chi_A$ is the map
\[
\chi_A : \GSp_4(\hat{\Z}) \to \hat{\Z}^\times \to (\Z/\abs{d_A}\Z)^\times \to \set{\pm 1}
\]
given by the multiplier map, followed by reduction modulo $\abs{d_A}$, followed by the kronecker symbol $\p{\frac{d_A}{\cdot}}$.

Assume that $A$ has the property that
\begin{equation} [\GSp_4(\hat{\Z}) : \rho_A(G_\Q)] = 2. \label{max-GA} \end{equation}
Then the inclusion in (\ref{epsilon-chi}) is an equality. As in the case for Serre curves \cite[Proposition 17]{MR4041152}, it then follows  that the image conductor for $A$ is given by
\[ m_A = \lcm(2,\abs{d_A}) = \begin{cases}
2\abs{\Delta_{\text{sf}}} & \Delta_{\text{sf}} \equiv 1 \pmod{4}, \\
4\abs{\Delta_{\text{sf}}} & \text{otherwise.}
\end{cases} \]
Thus if $A$ satisfies (\ref{max-GA}), the primes of bad reduction for $A$ and $C$ coincide and include $2$, and the discriminant $\Delta$ is square-free and congruent to $1$ modulo $4$, then  Theorem \ref{main-thrm} is sharp for $A$. The author is not aware of any such abelian surface in the literature, though an example satisfying (\ref{max-GA}) is given in \cite[Theorem 1.2]{MR4195609}.
\end{remark}

\begin{remark} The third condition in Theorem \ref{main-thrm} is rather unnatural. This assumption on $\mathcal{B}_A$ is used in the proof of Lemma \ref{ar-lem3}, and arises from the failure of a relevant lifting result in the case when $\ell = 2$ and $g = 2$. A careful analysis of $\GSp_4(\Z/8\Z)$ could perhaps lead to a refined condition (cf. \cite[pp. 13-14]{MR4190460}).
\end{remark}

\begin{remark} The constant $\mathcal{B}_A$ is constructed in view of Corollary \ref{NOS-cor}. Given this, it seems that one should be able write Theorem \ref{main-thrm} in terms of an arbitrary family of $G_K$-modules $\set{M[n]}_{n \geq 1}$ of $A(\overline{K})$ that satisfy the conclusion of Corollary \ref{NOS-cor}.
\end{remark}

\section{Notation and preliminaries}

\subsection{Symplectic groups \label{symp-1}} Let $R$ be a commutative ring with unity and let $M$ be a free $R$-module of rank $2g$. A map $\cyc{\cdot,\cdot}: M \oplus M \to R$ is called a \textit{symplectic form} on $M$ if it is bilinear, non-degenerate, and alternating. Given a symplectic form $\cyc{\cdot,\cdot}$ on $M$, the \textit{general symplectic group} and \textit{symplectic group} of $(M,\cyc{\cdot,\cdot})$ are
\begin{align*}
\GSp(M,\cyc{\cdot,\cdot}) &\coloneqq \set{\gamma \in \GL(M): \exists m(\gamma) \in R^\times \; \forall v,w \in M \; \cyc{\gamma v,\gamma w} = m(\gamma) \cyc{v,w} } \!, \\
\Sp(M,\cyc{\cdot,\cdot}) &\coloneqq \set{\gamma \in \GL(M): \forall v,w \in M \; \cyc{\gamma v,\gamma w} = \cyc{v,w}}\!.
\end{align*}

We may choose an $R$-basis for $M$ under which the symplectic form $\cyc{\cdot,\cdot}$ is represented by the block matrix
\[ \Omega_{2g} = 
\begin{pmatrix}
{0} & I_g \\ -I_g & {0} 
\end{pmatrix} \! , \]
where $I_g \in \Mat_{2g\times2g}(R)$ denotes the $g \times g$ identity matrix. Let $\mu: \GL(M) \overset{\sim}{\to} \GL_{2g}(R)$ be the isomorphism induced by our choice of basis. The images of $\GSp(M,\cyc{\cdot,\cdot})$ and $\Sp(M,\cyc{\cdot,\cdot})$ under $\mu$ are, respectively,
 \begin{align*}
 \GSp_{2g}(R) &\coloneqq \set{\gamma \in \GL_{2g}(R) : \exists m(\gamma) \in R^\times \text{ so that } \gamma^\intercal \Omega_{2g} \gamma =  m(\gamma) \Omega_{2g}}\!, \\
 \Sp_{2g}(R) &\coloneqq \set{\gamma \in \GL_{2g}(R) : \gamma^\intercal \Omega_{2g} \gamma =  \Omega_{2g}}\!.
 \end{align*}
The map $\mult: \GSp_{2g}(R) \twoheadrightarrow R^\times$ defined by $\gamma \mapsto m(\gamma)$ is a surjective homomorphism \cite[p. 50]{MR502254} and we see that
\[
\Sp_{2g}(R) = \ker\p{\GSp_{2g}(R) \xrightarrowdbl{\mult}  R^\times}\!.
\]
The orders of $\Sp_{2g}(R)$ and $\GSp_{2g}(R)$ are, in the important case of $R = \F_\ell$, given \cite[Theorem 3.1.2]{MR502254} by
\begin{equation} \abs{\Sp_{2g}(\F_\ell)} = \ell^{g^2} \prod_{i = 1}^{g} (\ell^{2i} - 1) \quad \text{and} \quad \abs{\GSp_{2g}(\F_\ell)} = (\ell-1)\ell^{g^2} \prod_{i = 1}^{g} (\ell^{2i} - 1). \label{symp-orders} \vspace{-6pt}
\end{equation}
\subsection{Notation} Throughout this paper, $p$ and $\ell$ denote prime numbers; $m$ and $n$ denote positive integers. 

Let $\hat{\Z}$ denote the ring of profinite integers and $\Z_\ell$ denote the ring of $\ell$-adic integers. The Chinese remainder theorem gives an isomorphism $\hat{\Z} \xrightarrow{\sim} \prod_\ell \Z_\ell$. The ring of $n$-adic integers $\Z_n$ and the ring of $(n)$-adic integers $\Z_{(n)}$ are, respectively, the quotients of $\hat{\Z}$ that correspond with $\Z_n \cong \prod_{\ell \mid n} \Z_\ell$ and $\Z_{(n)} \cong \prod_{\ell \nmid n} \Z_\ell$.

We see that $\hat{\Z} \cong \Z_n \times \Z_{(n)}$, and hence
\begin{equation} \GSp_{2g}(\hat{\Z}) \cong \GSp_{2g}(\Z_n) \times \GSp_{2g}(\Z_{(n)}). \label{gsp-zhat-iso} \end{equation}
Let $\rad(m) \coloneqq \prod_{\ell \mid m} \ell$ denote the radical of $m$. With (\ref{gsp-zhat-iso}) in mind, we define the following projection maps
\begin{align*}
\pi_n&: \GSp_{2g}(\hat{\Z}) \to \GSp_{2g}(\Z_n) \\
\pi_{(n)}&: \GSp_{2g}(\hat{\Z}) \to \GSp_{2g}(\Z_{(n)}) \\
\pi_{n^{\infty}, m}&: \GSp_{2g}(\Z_{n}) \to \GSp_{2g}(\Z/m\Z) && \text{(provided $\rad(m) \mid n$)}\\
\bar{\pi}_n&: \GSp_{2g}(\hat{\Z}) \to \GSp_{2g}(\Z/n\Z) \\
\pi_{n,m}&: \GSp_{2g}(\Z/n\Z) \to \GSp_{2g}(\Z/m\Z)  && \text{(provided $m \mid n$)}.
\end{align*}
For a closed subgroup $G \subseteq \GSp_{2g}(\hat{\Z})$, we employ the following notation
\[ G_n \coloneqq \pi_n(G), \qquad G_{(n)} \coloneqq \pi_{(n)}(G), \text{ and} \qquad G(n) \coloneqq \bar{\pi}_n(G). \]

Because Theorem \ref{main-thrm} is known \cite{MR4190460} for $g = 1$, in order to simplify our exposition, $g$ will always denote an integer that is at least two, unless otherwise stated. We shall often use the abbreviation $\ell_g$, which denotes
\begin{equation} \label{lg}
\ell_g \coloneqq 
\begin{cases}
3 & g = 2 \\
2 & g \geq 3
\end{cases}.
\end{equation} 

\subsection{Conductor\label{split-stable}} Let $G \subseteq \GSp_{2g}(\hat{\Z})$ be any open subgroup. Then $\set{\ker \bar{\pi}_m}_{m=1}^\infty$  is a neighborhood basis for the identity of $\GSp_{2g}(\hat{\Z})$. Hence, there exists an $m$ for which $\ker \bar{\pi}_m \subseteq G$. The \textit{conductor} of $G$ is
\begin{equation} m_G \coloneqq \min \set{m \in \N : \ker \bar{\pi}_m \subseteq G}\!. \label{tors-cond-gp} \end{equation}
It is sometimes helpful to understand the conductor  in the ways described in the following two lemmas.
\begin{lemma} \label{eq-defs} We have that $G = \bar{\pi}_m^{-1}(G(m))$ if and only if $\ker \bar{\pi}_m \subseteq G$. Consequently,
\[
m_G = \min \set{m \in \N : G = \bar{\pi}_{m}^{-1}(G(m))}\!.
\]
\end{lemma}
\begin{proof} We have $G \subseteq \bar{\pi}_m^{-1}(G(m))$, and both of these groups surject onto $G(m)$ via $\bar{\pi}_m$. Further, we see that
\[
\ker\p{\bar{\pi}_m^{-1}(G(m)) \xrightarrowdbl{\bar{\pi}_m} G(m)} =  \ker \bar{\pi}_m \quad \text{and} \quad
\ker\p{G \xrightarrowdbl{\bar{\pi}_m} G(m)} = G \cap \ker \bar{\pi}_m.
\]
Thus, $G = \bar{\pi}_m^{-1}(G(m))$ if and only if $\ker \bar{\pi}_m = G \cap \ker \bar{\pi}_m$, which happens if and only if $\ker \bar{\pi}_m \subseteq G$.
\end{proof}

For our next lemma, we give some terminology (see, \cite[I \S 1.1]{MR0568299}). We say that $m$ \textit{splits} $G$ if
\[ (\pi_m \times \pi_{(m)})(G) = G_m \times \GSp_{2g}(\Z_{(m)}). \]
We say that $m$ is \textit{stable} for $G$ if
\[ G_m = \pi_{m^\infty,m}^{-1}(G(m)). \]

\begin{lemma} \label{eq-defs2} We have that $G = \bar{\pi}_m^{-1}(G(m))$ if and only if $m$ splits and is stable for $G$. Consequently,
\[
m_G = \min \set{m \in \N : m \text{ splits and is stable for } G}\!.
\]
\end{lemma}
\begin{proof} The map $\pi_m \times \pi_{(m)}: \GSp_{2g}(\hat{\Z}) \to \GSp_{2g}(\Z_m) \times \GSp_{2g}(\Z_{(m)})$ is an isomorphism. We see that
\[ (\pi_m \times \pi_{(m)})(\bar{\pi}_m^{-1}(G(m))) =  \pi_{m^\infty,m}^{-1}(G(m)) \times \GSp_{2g}(\Z_{(m)}). \]
Thus $G = \bar{\pi}_m^{-1}(G(m))$ if and only if $m$ splits and is stable for $G$. The conclusion follows from Lemma \ref{eq-defs}.
\end{proof}

\subsection{Galois representations} Let $A$ be a principally polarized abelian variety of dimension $g$ over a number field $K$. Let $T(A) \coloneqq \varprojlim A[m]$ be the adelic Tate module of $A$. Recall that $T(A)$ is a free $\hat{\Z}$-module of rank $2g$. The Weil pairing and a choice of principal polarization on $A$ yield a symplectic form $\cyc{\cdot,\cdot}: T(A) \oplus T(A) \to \hat{\Z}^\times$. The continuous action of $G_K$ on $T(A)$ is compatible with this symplectic form and hence induces a representation $G_K \to \GSp(T(A), \cyc{\cdot,\cdot})$. With a choice of basis, we obtain the continuous homomorphism of profinite groups
\[\rho_A : G_K \to \GSp_{2g}(\hat{\Z}) \]
known as the \textit{adelic Galois representation of $A$}. The \textit{Galois image} of $A$ is the subgroup $G \coloneqq \rho_A(G_K)$ of $\GSp_{2g}(\hat{\Z})$. If $G$ is open in $\GSp_{2g}(\hat{\Z})$, the \textit{image conductor} of $A$ is defined to be the conductor of $G$ as in (\ref{tors-cond-gp}).

\begin{remark} \label{galrep-rmks}  Below are three key observations relating to the Galois image $G$ of $A$.

\begin{enumerate}
\item We see that $G$ is a closed subgroup of the profinite group $\GSp_{2g}(\hat{\Z})$. A consequence is that $G$ is an open subgroup of $\GSp_{2g}(\hat{\Z})$ if and only if the group index $[\GSp_{2g}(\hat{\Z}) : G]$ is finite.
\item  For a subset $S \subseteq A(\overline{K})$, let $K(S)$ be the extension of $K$ obtained by adjoining to $K$ the coordinates of the points in $S$. Let $A[n]$ be the $n$-torsion subgroup of $A(\overline{K})$, $A[n^\infty] \coloneqq \bigcup_{k=0}^\infty A[n^k]$, and $A_\tors \coloneqq \cup_{n=1}^\infty A[n]$. We have
\begin{align*}
G &\cong \Gal(K(A_{\tors})/K) \\
G_n &\cong \Gal(K(A[n^\infty])/ K) \\
G(n) &\cong \Gal(K(A[n]) / K).
\end{align*}
Further, let $A_{\tors,(n)} \coloneqq \bigcup_{\gcd(m,n)=1} A[m]$. We have that 
\[ G_{(n)} \cong \Gal(K(A_{\tors,(n)})/K).\]
\item Let $\mu_n$ be the group of $n$th roots of unity in $\overline{K}$. Let $\mu_{\ell^\infty} \coloneqq \bigcup_k \mu(\ell^k)$ and $\mu_\infty \coloneqq \bigcup_n \mu(n)$. The composition $\mult \circ \rho_A: G_K \to \GSp_{2g}(\hat{\Z}) \to \hat{\Z}^\times$ is the cyclotomic character of $K$. Thus,
\begin{align*}
\mult(G) &\cong \Gal(K(\mu_\infty)/K) \\
\mult(G_\ell) &\cong \Gal(K(\mu_{\ell^\infty})/K) \\
\mult(G(n)) &\cong \Gal(K(\mu_n) / K).
\end{align*}
\end{enumerate}
\end{remark}

We now give a generalization of a variant of \cite[IV-18 Lemma (2)]{MR1484415}.

\begin{lemma} \label{ar-lem2} As before, let $G \coloneqq \rho_A(G_K)$. If $\ell$ is such that $\Sp_{2g}(\Z_\ell) \subseteq G_\ell$, then $G_\ell = \GSp_{2g}(\Z_\ell)$ if and only if $K \cap \Q(\mu_{\ell^\infty}) = \Q$. In particular, if $\Sp_{2g}(\Z_\ell) \subseteq G_\ell \neq \GSp_{2g}(\Z_\ell)$, then $K/\Q$ is ramified at $\ell$.
\end{lemma} 
\begin{proof} Since $\Sp_{2g}(\Z_\ell) \subseteq G_\ell$, both $\mult: \GSp_{2g}(\Z_\ell) \twoheadrightarrow \Z_\ell^\times$ and the restriction $\mult|_{G_\ell}: G_\ell \to \Z_\ell^\times$ have kernel $\Sp_{2g}(\Z_\ell)$. Therefore, $G_\ell = \GSp_{2g}(\Z_\ell)$ if and only if $\mult(G_\ell) = \Z_\ell^\times$. By Remark \ref{galrep-rmks}(3) and Galois theory,
\[ \mult(G_\ell) \cong \Gal(K(\mu_{\ell^\infty})/K) \cong \Gal(\Q(\mu_{\ell^\infty})/(K \cap \Q(\mu_{\ell^\infty}))) \subseteq \Gal(\Q(\mu_{\ell^\infty})/\Q) \cong \Z_\ell^\times. \]
It follows that $\mult(G_\ell) = \Z_\ell^\times$ if and only if the extension $K \cap \Q(\mu_{\ell^\infty}) / \Q$ is nontrivial. 

Now assume $\Sp_{2g}(\Z_\ell) \subseteq  G_\ell \neq \GSp_{2g}(\Z_\ell)$. By the above, the extension $K \cap \Q(\mu_{\ell^\infty}) / \Q$ is nontrivial. Thus, this extension is ramified at $\ell$ as it is a sub-extension of $\Q(\mu_{\ell^\infty})/\Q$, which is well-known to be totally ramified at $\ell$. Thus, $K/\Q$ is ramified at $\ell$ because it has a ramified sub-extension.
\end{proof}

\subsection{Fiber product}  Let $G_1, G_2,$ and $Q$ be groups. Let $\psi_1: G_1 \twoheadrightarrow Q$ and $\psi_2: G_2 \twoheadrightarrow Q$ be surjective homomorphisms. The \textit{fiber product} of $G_1$ and $G_2$ over $(\psi_1,\psi_2)$ is the group
\[
G_1 \times_{(\psi_1,\psi_2)} G_2 \coloneqq \set{(g_1,g_2) \in G_1 \times G_2 : \psi_1(g_1) = \psi_2(g_2)}\!.
\]
Observe that $G_1 \times_{(\psi_1,\psi_2)} G_2 \subseteq G_1 \times G_2$ is a subgroup that surjects onto both $G_1$ and $G_2$ via the relevant projection maps. Writing $\psi = (\psi_1,\psi_2)$, we say that a fiber product $G_1 \times_{\psi} G_2$ is \textit{trivial} if $G_1 \times_{\psi} G_2 = G_1 \times G_2$.

Let $L_1/K$ and $L_2/K$ be Galois extensions, both contained in $\overline{K}$. The \textit{entanglement field} of $L_1$ and $L_2$ is the intersection $L_1 \cap L_2$. The \textit{compositum} of $L_1$ and $L_2$, denoted $L_1 L_2$, is the smallest (by inclusion) subfield of $\overline{K}$ containing both $L_1$ and $L_2$. The Galois group of $L_1L_2/K$ may be described using the fiber product.

\begin{lemma} Let $L_1 / K$ and $L_2 / K$ be Galois extensions, contained in $\overline{K}$. Then $L_1L_2/K$ is Galois and
\[ \Gal(L_1L_2/K) \cong \Gal(L_1/K) \times_{(\psi_1,\psi_2)} \Gal(L_2/K), \]
where each $\psi_i : \Gal(L_i/K) \twoheadrightarrow \Gal(L_1 \cap L_2  / K)$ is the canonical restriction homomorphism. \label{goursat-field}
\end{lemma}
\begin{proof} See \cite[Theorem VI 1.14]{MR1878556}.
\end{proof}

\section{Symplectic groups}

In \S\ref{symp-1}, we introduced the symplectic groups $\GSp_{2g}(R)$ and $\Sp_{2g}(R)$. In this section, we  derive some useful properties of these groups when $R = \F_\ell$ and $R = \Z_\ell$.

\subsection{Normal subgroups} The objective of this subsection is to understand the normal subgroups of $\GSp_{2g}(\F_\ell)$ for $\ell \geq \ell_g$, where $\ell_g$ is as in (\ref{lg}). We begin by considering the projective symplectic groups.

The  center of $\GSp_{2g}(\F_\ell)$ is the scalar subgroup $\Lambda_{2g}(\F_\ell)$ of $\GL_{2g}(\F_\ell)$ \cite[4.2.5(5)]{MR502254}. Let $\pi$ be the projection
\[ \pi: \GSp_{2g}(\F_\ell) \twoheadrightarrow \GSp_{2g}(\F_\ell) / \Lambda_{2g}(\F_\ell). \]
The \textit{projective general symplectic group} $\PGSp_{2g}(\F_\ell)$ and \textit{projective symplectic group} $\PSp_{2g}(\F_\ell)$ are the images of $\GSp_{2g}(\F_\ell)$ and $\Sp_{2g}(\F_\ell)$ under $\pi$, respectively. We give some useful properties of these groups below. Here and later, we let $[\cdot,\cdot]$ denote a commutator and write $G'$ to denote the commutator subgroup of a group $G$.

\begin{lemma} \label{symplectic-facts} Assume $\ell \geq \ell_g$. Each of the following statements hold.
\begin{enumerate}
\item The center of $\PGSp_{2g}(\F_\ell)$ is trivial;
\item $\GSp_{2g}(\F_\ell)' = \Sp_{2g}(\F_\ell)' = \Sp_{2g}(\F_\ell)$;
\item $\PGSp_{2g}(\F_\ell)' =  \PSp_{2g}(\F_\ell)$; and
\item $\PSp_{2g}(\F_\ell)$ is simple.
\end{enumerate}
\end{lemma}
\begin{proof} Statements (1), (2), and (4) are found in \cite[4.2.5(2), 3.3.6, 3.4.1]{MR502254}. For (3), we apply (2) to see that 
\[ \PGSp_{2g}(\F_\ell)' = \p{\pi(\GSp_{2g}(\F_\ell))}' = \pi(\GSp_{2g}(\F_\ell)') = \pi(\Sp_{2g}(\F_\ell)) = \PSp_{2g}(\F_\ell). \qedhere \]
\end{proof}

Using the properties of Lemma \ref{symplectic-facts}, we now determine the normal subgroups of $\PGSp_{2g}(\F_\ell)$. Our target lemma regarding the normal subgroups of $\GSp_{2g}(\F_\ell)$ then follows. We make the abbreviation $\Lambda_{2g} \coloneqq \Lambda_{2g}(\F_\ell)$.

\begin{lemma} \label{norm-pgsp} Assume that $\ell \geq \ell_g$. If $N \normal \PGSp_{2g}(\F_\ell)$, then either $N = \set{\Lambda_{2g}}$ or $\PSp_{2g}(\F_\ell) \subseteq N$.

\end{lemma}
\begin{proof} Assume that $N \normal \PGSp_{2g}(\F_\ell)$ is nontrivial. Since the center of $\PGSp_{2g}(\F_\ell)$ is trivial, we have
\[ \set{\Lambda_{2g}} \subsetneq [\PGSp_{2g}(\F_\ell), N] \subseteq N \cap \PGSp_{2g}(\F_\ell)' = N \cap \PSp_{2g}(\F_\ell) \normal \PSp_{2g}(\F_\ell). \]
By the simplicity of $\PSp_{2g}(\F_\ell)$, this implies that $N \cap \PSp_{2g}(\F_\ell) = \PSp_{2g}(\F_\ell)$. Thus, $\PSp_{2g}(\F_\ell) \subseteq N$.
\end{proof}

\begin{lemma} \label{norm-gsp} Assume that $\ell \geq \ell_g$. If $N \normal \GSp_{2g}(\F_\ell)$, then either $N \subseteq \Lambda_{2g}$ or $\Sp_{2g}(\F_\ell) \subseteq N$.
\end{lemma}
\begin{proof} Assume that $N \not\subseteq \Lambda_{2g}$. Then $\pi(N) \normal \PGSp_{2g}(\F_\ell)$ is nontrivial. So, by Lemma \ref{norm-pgsp}, $\PSp_{2g}(\F_\ell) \subseteq \pi(N)$ and hence $\Sp_{2g}(\F_\ell) \Lambda_{2g} \subseteq N \Lambda_{2g}$. By taking commutators, we now see that
\[    N \supseteq N' = (N \Lambda_{2g})' \supseteq (\Sp_{2g}(\F_\ell) \Lambda_{2g})' = (\Sp_{2g}(\F_\ell))'= \Sp_{2g}(\F_\ell).   \qedhere \]
\end{proof}

\subsection{Index bound\label{indx-bd}} Here we use Lemma \ref{norm-gsp} and a standard lemma from group theory to obtain a lower bound on the index of each subgroup of $\GSp_{2g}(\F_\ell)$ that does not contain $\Sp_{2g}(\F_\ell)$. We write $n!$ to denote the factorial of $n$.

\begin{lemma} \label{norm-core} Let $G$ be a finite group and $H \subseteq G$ a subgroup. The \textit{normal core} of $H$ in $G$, denoted $H_G$, is the largest (by inclusion) subgroup of $H$ that is normal in $G$. One has that $[G: H_G]$ divides $[G:H]!$.
\end{lemma}
\begin{proof} See \cite[1.6.9]{MR1357169}.
\end{proof}

\begin{lemma} \label{gp-lem} Let $G \subseteq \GSp_{2g}(\F_\ell)$ be a subgroup. If $\Sp_{2g}(\F_\ell) \not\subseteq G$, then 
\[ \left[ \GSp_{2g}(\F_\ell) : G \right] \geq \ell. \]
\end{lemma}
\begin{proof} The result is clear for $\ell = 2$, so we assume that $\ell \geq 3$. Let $N$ be the normal core of $G$ in $\GSp_{2g}(\F_\ell)$. Then $N \normal \GSp_{2g}(\F_\ell)$ and $\Sp_{2g}(\F_\ell) \not\subseteq N$, so $N \subseteq \Lambda_{2g}(\F_\ell)$, by Lemma \ref{norm-gsp}. Now, by (\ref{symp-orders}) and Lemma \ref{norm-core},
\[\ell \divides \abs{\PGSp_{2g}(\F_\ell)} \divides [\GSp_{2g}(\F_\ell):N] \divides [\GSp_{2g}(\F_\ell): G]!. \qedhere \]
\end{proof}

\subsection{Subgroup lifting}

We state a lifting lemma for $\Sp_{2g}(\Z_\ell)$ that extends \cite[IV-23 Lemma 3]{MR1484415}. Then, we give two corollaries and state a lifting lemma $\GSp_{2g}(\Z_\ell)$. As before, we shall assume that $g \geq 2$.

\begin{proposition} \label{sp-lift} Let $H_\ell \subseteq \Sp_{2g}(\Z_\ell)$ be a closed subgroup. If $H(\ell) = \Sp_{2g}(\Z/\ell\Z)$, then $H_\ell = \Sp_{2g}(\Z_\ell)$.
\end{proposition}
\begin{proof} See \cite[Theorem 1]{MR3667841}.
\end{proof}

For a subgroup $H \subseteq G_\ell$, we let $\overline{H}$ denote the topological closure of $H$ in $G_\ell$.

\begin{corollary} \label{sp-gsp-lift} Assume that $\ell \geq \ell_g$ and let $G_\ell \subseteq \GSp_{2g}(\Z_\ell)$ be a closed subgroup. If $\Sp_{2g}(\Z/\ell\Z) \subseteq G(\ell)$, then  $\Sp_{2g}(\Z_\ell) \subseteq G_\ell$.
\end{corollary}
\begin{proof}   We have that $\overline{(G_\ell)'} \subseteq \Sp_{2g}(\Z_\ell)$ is a closed subgroup. Further, as $G_\ell$ surjects onto $G(\ell)$, we have 
\[ (G_\ell)'(\ell) = (G(\ell))' = \Sp_{2g}(\Z/\ell\Z), \]
 by Lemma \ref{symplectic-facts}(2).  Thus, $\overline{(G_\ell)'}(\ell)  = \Sp_{2g}(\Z/\ell\Z)$. So, by Proposition \ref{sp-lift}, $ G_\ell \supseteq \overline{(G_\ell)'} = \Sp_{2g}(\Z_\ell)$. 
\end{proof}

\begin{corollary} \label{normal-mult-gsp} Assume that $\ell \geq 3$ and let $N_\ell \normal \GSp_{2g}(\Z_\ell)$ be a closed normal subgroup. If $\mult(N_\ell) = \Z_\ell^\times$, then $N_\ell = \GSp_{2g}(\Z_\ell)$.
\end{corollary}
\begin{proof} Since $\mult(N(\ell)) = (\Z/\ell\Z)^\times$ and $\mult(\Lambda_{2g}(\F_\ell)) = (\Z/\ell\Z)^{\times 2}$, we have $N(\ell) \not\subseteq \Lambda_{2g}(\F_\ell)$. Hence, $\Sp_{2g}(\Z/\ell\Z) \subseteq N(\ell)$  by Lemma \ref{norm-gsp}. Thus, by Corollary \ref{sp-gsp-lift}, we find that $\Sp_{2g}(\Z_\ell) \subseteq N_\ell$. As both $\GSp_{2g}(\Z_\ell)$ and $N_\ell$ surject onto $\Z_\ell^\times$, via $\mult$, with kernel $\Sp_{2g}(\Z_\ell)$, we conclude  that  $N_\ell = \GSp_{2g}(\Z_\ell)$.
\end{proof}

We now state a lifting lemma for $\GSp_{2g}(\Z_\ell)$. Let $\alpha_\ell$ denote the quantity
\begin{equation} \label{alphal}
\alpha_\ell \coloneqq \begin{cases} 2 & \text{if } \ell = 2 \\  1 & \text{if } \ell \geq 3 \end{cases}.
\end{equation}
\begin{lemma} \label{GSp-lift} Let $G_\ell \subseteq \GSp_{2g}(\Z_\ell)$ be a closed subgroup. We have that if $G(\ell^{\alpha_\ell+1}) = \GSp_{2g}(\Z/\ell^{\alpha_\ell+1}\Z)$, then $G_\ell = \GSp_{2g}(\Z_\ell)$.
\end{lemma}
\begin{proof} See \cite[Remark 3.2]{MR4190460}, the proof of which generalizes directly to arbitrary $g$, mutatis mutandis. 
\end{proof}

\section{Proof of Theorem \ref{main-thrm}, assuming two propositions \label{props-imply-thrm}}

We begin by stating two propositions, which we shall prove in \S\ref{prop-1-sec} and \S\ref{prop-2-sec}.  The first proposition is purely group-theoretic, whereas the second depends on the arithmetic of the abelian variety $A$. Due to group-theoretic differences relating to the prime 2 (visible in the statement of Lemma \ref{GSp-lift}), we employ a variant of the radical function. This modified radical is denoted $\rad'$ and is defined by
\[
\rad'(n) \coloneqq \begin{cases}
2\rad(n) & 4 \mid n \\
\rad(n) & \text{otherwise},
\end{cases}
\]
where $\rad(n) = \prod_{\ell \mid n} \ell$ is the usual radical of $n$. Our main propositions are as follows.

\begin{proposition} \label{prop-1} Let $g$ be an integer, $G \subseteq \GSp_{2g}(\hat{\Z})$ be an open subgroup, and $m_G$ be as in (\ref{tors-cond-gp}). Then
\[ \frac{m_G}{\rad'(m_G)} \text{ divides } \left[ \pi_{m_G, \rad'(m_G)}^{-1}(G(\rad'(m_G))) : G(m_G) \right]\!. \]
\end{proposition}

\begin{proposition}\label{prop-2} Let $g \geq 2$ be an integer and let $A$ be as in the statement of Theorem \ref{main-thrm}. Then
\[ \rad'(m_{A}) \leq 2 \cdot  \mathcal{B}_A \cdot \left[ \GSp_{2g}(\Z/\rad'(m_A)\Z) : G(\rad'(m_A)) \right]\!, \]
where $G$ is the Galois image of $A$, $m_A$ is the image conductor of $A$, and $\mathcal{B}_{A}$ is as in Theorem \ref{main-thrm}. 
\end{proposition}

We now prove Theorem \ref{main-thrm}, assuming Proposition \ref{prop-1} and Proposition \ref{prop-2}.
\begin{proof}[Proof of Theorem \ref{main-thrm}] Write $G \coloneqq \rho_{A}(G_K)$ and  $r' \coloneqq \rad'(m_A)$. Using Lemma \ref{eq-defs} initially, we see
\begin{align*}
[\GSp_{2g}(\hat{\Z}) : G] &= [\GSp_{2g}(\hat{\Z}) : \bar{\pi}_{m_A}^{-1}(G(m_{A}))] \\
&= [\GSp_{2g}(\Z/m_{A}\Z) : G(m_{A})] \\
&= [\GSp_{2g}(\Z/m_{A}\Z) : \pi_{m_{A},r'}^{-1}(G(r'))] [\pi_{m_{A},r'}^{-1}(G(r')) : G(m_{A})] \\
&= [\GSp_{2g}(\Z/r'\Z) : G(r')] [\pi_{m_{A},r'}^{-1}(G(r')) : G(m_{A})].
\end{align*}
With the above in mind, applying Proposition \ref{prop-1} and Proposition \ref{prop-2} now yields
\begin{align*}
m_{A} &= r' \cdot \frac{m_{A}}{r'}  \\
&\leq 2 \cdot \mathcal{B}_A \cdot  \left[ \GSp_{2g}(\Z/r'\Z) : G(r') \right] \cdot [ \pi_{m_A, r'}^{-1}(G(r')) : G(m_A) ]  \\ 
&= 2\cdot  \mathcal{B}_A \cdot[\GSp_{2g}(\hat{\Z}):G]. \qedhere
\end{align*}
\end{proof}

\section{Proof of Proposition \ref{prop-1} \label{prop-1-sec}} 

For the case of $g = 1$, a proof of Proposition \ref{prop-1} is given in \cite[Proposition 1.6]{MR4190460}. This purely group theoretic proof immediately generalizes, mutatis mutandis, to prove Proposition \ref{prop-1} for arbitrary $g$. For this reason, in this section we shall explain the structure of the proof and refer the reader to \cite{MR4190460} for the details.

Let $G \subseteq \GSp_{2g}(\hat{\Z})$ be any open subgroup and write $m_G \eqqcolon \prod_{\ell \mid m_G} \ell^{\beta_\ell}$ for the prime factorization of its conductor. For each $k$, write $N_{\ell^k} \coloneqq \ker(\pi_{\ell^{k+1},\ell^k})$. Using a lifting lemma \cite[Lemma 3.1]{MR4190460}, we may describe \cite[Corollary 3.5]{MR4190460} each $\beta_\ell$ as
\[
\beta_\ell = \min\set{
  \beta \geq 0 
  :
  \forall k \in [\beta,\max\set{\beta,\alpha_\ell}] \cap \Z, \,
  N_{\ell^k} \times \{1_{(\ell)}\}
  \subseteq
  (\pi_{\ell^\infty,\ell^{k+1}} \times \pi_{(\ell)})(G)
  } \!,
\]
where $\alpha_\ell$ is defined in (\ref{alphal}) and $1_{(\ell)}$ denotes the identity of $\GSp_{2g}(\Z_{(\ell)})$. As a corollary, it follows \cite[Lemma 3.8]{MR4190460} that if $d$ is a positive integer that satisfies the divisibility condition  $\rad'(m_G) \mid d \mid d\ell \mid m_G $, then 
\begin{equation} \ell \divides [\pi_{\ell d,d}^{-1}(G(d)): G(\ell d)]. \label{J-L3.7}
\end{equation}

Write $r'\coloneqq \rad'(m_G)$. Let $\ell$ be a prime dividing $\frac{m_G}{r'}$. Let $\beta_\ell$ and $r_\ell$ be such that $\ell^{\beta_\ell} \mid\mid m_G$ and $\ell^{r_\ell} \mid\mid r'$, respectively. Applying (\ref{J-L3.7}) with $d = \ell^k r'$ for each integer $k$ such that $0 \leq k < \beta_\ell - r_\ell$, we obtain that
\begin{align*}
\ell^{\beta_\ell - r_\ell} &\divides \prod_{0 \leq k < \beta_\ell - r_\ell} [\pi_{\ell^{k+1}r' ,\ell^kr'}^{-1}(G(\ell^kr')): G(\ell^{k+1} r')]  \\
&\divides [\pi^{-1}_{\ell^{\beta_\ell - r_\ell}r',r'}(G(r')) : G(\ell^{\beta_\ell - r_\ell}r')] \\
&\divides [\pi^{-1}_{m_G,r'}(G(r')): G(m_G)].
\end{align*}
Since the above holds for each prime $\ell$ dividing $\frac{m_G}{\rad'(m_G)}$, it follows that
\[ \frac{m_G}{\rad'(m_G)} = \prod_{\ell \mid \frac{m_G}{r'}} \ell^{\beta_\ell - r_\ell}  \divides  [\pi^{-1}_{m_G,\rad'(m_G)}(G(\rad'(m_G))): G(m_G)]. \]

\section{Constraints on prime divisors of the image conductor} \label{constraints}

Let $A$ be as in the statement of Theorem \ref{main-thrm}. We give constraints on the primes that divide the image conductor of $A$. To do so, we employ a variant of the N\'eron-Ogg-Shafarevich criterion for abelian varieties.

\begin{theorem}[Serre-Tate, 1968 {\cite{MR236190}}] \label{NOS} Let $A$ be an abelian variety over a number field $K$. Let $\mathcal{L} \subseteq \mathcal{O}_K$ be a prime ideal of $K$, lying over a rational prime $\ell$. The following are equivalent:
\begin{enumerate}
  \item $A$ has good reduction at $\mathcal{L}$;
  \item For each positive integer $m$ that is not divisible by $\ell$, the prime $\mathcal{L}$ is unramified in $K(A[m])/K$; and
  \item The prime $\mathcal{L}$ is unramified in $K(A_{\tors,(\ell)})/K$, where $A_{\tors,(\ell)}$ is defined in Remark \ref{galrep-rmks}(2).
\end{enumerate}
\end{theorem}

Recall that the constant $\mathcal{B}_A$ is defined in the statement of Theorem \ref{main-thrm}.

\begin{corollary} Assume that $\ell \geq \ell_g$. Then $\ell$ divides $\mathcal{B}_A$ if and only if $K(A_{\tors,(\ell)})/\Q$ is ramified at $\ell$. \label{NOS-cor}
\end{corollary}
\begin{proof} Since $\ell \geq \ell_g$, we have that $\ell$ divides $\mathcal{B}_A$ if and only if $K/\Q$ is ramified at $\ell$ or $A$ has bad reduction at some prime ideal of $\mathcal{O}_K$ that lies over $\ell$. By Theorem \ref{NOS}, the second condition is equivalent to the condition that $K(A_{\tors,(\ell)})/K$ is ramified at some prime ideal of $\mathcal{O}_K$ that lies over $\ell$. 
\end{proof}

Recall the notation of Remark \ref{galrep-rmks}(2) and that $G$ denotes the Galois image of $A$. The next lemma is key. It uses our understanding of $\mathcal{B}_A$ from Corollary \ref{NOS-cor}  to give a constraint on odd primes $\ell$ that divide $m_A$ for which $G_\ell = \GSp_{2g}(\Z_\ell)$.

\begin{lemma} \label{ar-lem1} Let $\ell$ be an odd prime that divides  $m_A$. If $G_\ell = \GSp_{2g}(\Z_\ell)$, then  $\ell$ divides $\mathcal{B}_A$.
\end{lemma}
\begin{proof} We see that $m_A \neq \ell$ for otherwise, by Lemma \ref{eq-defs}, we would have that $G = \pi_{\ell}^{-1}(\GSp_{2g}(\Z_\ell)) = \GSp_{2g}(\hat{\Z})$ and hence $m_A = 1$ is not divisible by $\ell$. Thus, as $\ell$ is stable for $G$, it follows from Lemma \ref{eq-defs2} that $\ell$ does not split $G$. Let $F$ be the entanglement field  $F \coloneqq K(A[\ell^\infty])  \cap K(A_{\tors,(\ell)})$. Then, by Lemma \ref{goursat-field}, $G$ may be expressed as the \textit{nontrivial} fiber product
\[
G \cong \Gal(K(A[\ell^\infty])K(A_{\tors,(\ell)})/K) \cong \GSp_{2g}(\Z_\ell) \times_{(\psi_\ell,\psi_{(\ell)})} \Gal(K(A_{\tors,(\ell)})/K), \label{G-fiber}
\]
where $\psi_\ell$ and $\psi_{(\ell)}$ are, upon making the identifications of  Remark \ref{galrep-rmks}(3), the restriction homomorphisms
\[\psi_{\ell}: \GSp_{2g}(\Z_\ell)  \twoheadrightarrow  \Gal(F/K) \quad \text{and} \quad
\psi_{(\ell)}: \Gal(K(A_{\tors,(\ell)})/K) \twoheadrightarrow  \Gal(F/K). \]
As the fiber product is nontrivial, in particular $\Gal(F/K)$ is nontrivial. Consider the following field diagram.
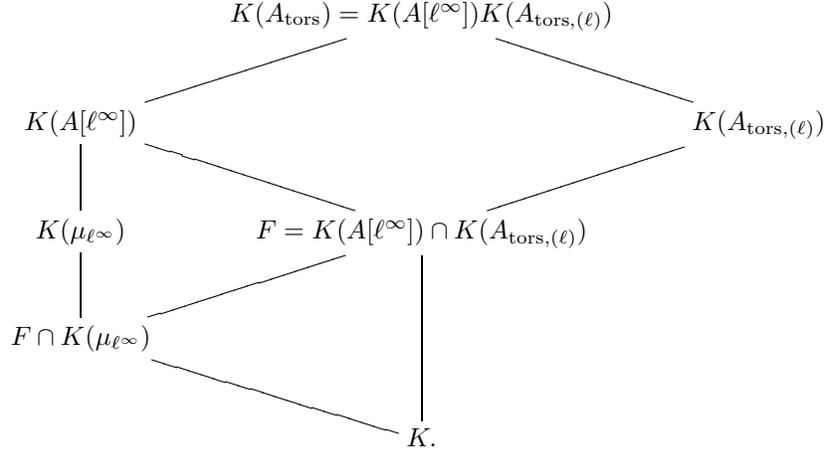
\[
\xymatrix{
  & K(A_{\tors}) = K(A[\ell^\infty])K(A_{\tors,(\ell)}) \ar@{-}[dl] \ar@{-}[dr] & \\
  K(A[\ell^\infty])\ar@{-}[d]   \ar@{-}[dr] & & K(A_{\tors,(\ell)}) \ar@{-}[dl] \\
  K(\mu_{\ell^\infty})\ar@{-}[d]  & F = K(A[\ell^\infty])  \cap K(A_{\tors,(\ell)})\ar@{-}[dl]  \ar@{-}[dd] & \\
  F \cap K(\mu_{\ell^\infty}) \ar@{-}[dr]&  & \\
  & K. & }
\]

If $K/\Q$ is ramified at $\ell$, then $\ell$ divides $\mathcal{B}_A$, so we are done. As such, we assume $K/\Q$ is unramified at $\ell$. Note that then $K(\mu_{\ell^\infty})/K$ is totally ramified at each prime ideal of $\mathcal{O}_K$ that lies over $\ell$. To show that $\ell$ divides $\mathcal{B}_A$, it suffices by Corollary \ref{NOS-cor} to show that $K(A_{\tors,(\ell)})/K$ is ramified at some prime ideal of $\mathcal{O}_K$ that lies over $\ell$. Hence, it suffices merely to show that the extension $F \cap K(\mu_{\ell^\infty}) / K$ is nontrivial.

Because $\psi_\ell$ is a surjective group homomorphism with nontrivial image, its kernel $\ker(\psi_\ell) \triangleleft \GSp_{2g}(\Z_\ell)$ is proper. Thus, by Corollary \ref{normal-mult-gsp}, we have that $\mult(\ker\psi_\ell)$ is a proper subgroup of $\Z_\ell^\times$. So we see,
\[ \mult(\cyc{\ker\psi_\ell,\Sp_{2g}(\Z_\ell)}) = \mult(\ker \psi_\ell) \subsetneq \Z_\ell^\times. \]
Hence $\cyc{\ker\psi_\ell,\Sp_{2g}(\Z_\ell)}  \triangleleft \GSp_{2g}(\Z_\ell)$ is a proper subgroup. Thus, by Galois theory,
\begin{align*} 
F \cap K(\mu_{\ell^\infty}) &= K(A[\ell^\infty])^{\ker\psi_\ell} \cap K(A[\ell^\infty])^{\Sp_{2g}(\Z_\ell)} \\
&= K(A[\ell^\infty])^{\cyc{\ker\psi_\ell,\Sp_{2g}(\Z_\ell)}} \\
&\supsetneq  K(A[\ell^\infty])^{\GSp_{2g}(\Z_\ell)} \\
&= K.
\end{align*}
We see that the extension $F \cap K(\mu_{\ell^\infty}) / K$  is nontrivial, and hence $\ell$ divides $\mathcal{B}_A$.
\end{proof}

Following Lemma \ref{ar-lem1}, which considers a odd prime $\ell$, our next lemma offers a constraint when $\ell = 2$ divides $m_A$.

\begin{lemma} \label{ar-lem3} Assume that $2$ divides $m_A$. Write $r' \coloneqq \rad'(m_A)$ and $s \coloneqq \frac{1}{2}r'$. We have
\begin{align}
(\bar{\pi}_2 \times \bar{\pi}_s)(G) \neq \GSp_{2g}(\Z/2\Z) \times G(s) &\implies 2 \leq  \left[ \pi_{r',s}^{-1}\p{G(s)}: G(r') \right]\!, \label{A2:1} \\
(\bar{\pi}_2 \times \bar{\pi}_s)(G) = \GSp_{2g}(\Z/2\Z) \times G(s) &\implies  2 \text{ divides } \mathcal{B}_A. \label{A2:2}
\end{align}
\end{lemma}
\begin{proof} Assume first that the hypothesis of (\ref{A2:1}) holds. Then
\[
(\pi_{r',2} \times \pi_{r',s})( G(r')) \neq \GSp_{2g}(\Z/2\Z) \times G(s) = (\pi_{r',2} \times \pi_{r',s})( \pi_{r',s}^{-1}\p{G(s)}).
\]
Thus $G(r')$ is a proper subgroup of $\pi_{r',s}^{-1}\p{G(s)}$, so the conclusion of (\ref{A2:1}) follows.

Now assume that the hypothesis of (\ref{A2:2}) holds. If $g = 2$, then 2 divides $\mathcal{B}_A$ by definition. As such, we assume $g \geq 3$. By hypothesis, $G(2) = \GSp_{2g}(\Z/2\Z)$, so $\Sp_{2g}(\Z_2) \subseteq G_2$ by Corollary \ref{sp-gsp-lift}. If $G_2 \neq \GSp_{2g}(\Z_2)$, then by Lemma \ref{ar-lem2}, the prime $2$ is ramified in $K/\Q$, so $2$ divides $\mathcal{B}_A$. Assume, therefore, that $G_2 = \GSp_{2g}(\Z_2)$.

We have that $2$ properly divides $m_A$ by Lemma \ref{eq-defs}. Thus, it follows from Lemma \ref{eq-defs2} that
\[ G_{r'} \cong (\pi_{2} \times \pi_{s})(G) = \GSp_{2g}(\Z_2) \times_\psi G_s
\]
is a nontrivial fiber product. Observe that each nontrivial finite quotient of $\GSp_{2g}(\Z_2)$ has even order whereas each nontrivial finite quotient of $\ker(\pi_{s^\infty,s})$ has odd order. For this reason, the fiber product
\[ (\pi_{2} \times \bar{\pi}_{s})(G) = \GSp_{2g}(\Z_2) \times_\psi G(s) \]
is nontrivial as well. Making the identifications of Remark \ref{galrep-rmks}(3), we conclude that the entanglement field $F \coloneqq K(A[2^\infty]) \cap K(A[s])$ is a nontrivial extension of $K$.

Consider the Galois group $H \coloneqq \Gal(K(A_\tors)/F)$. As $F/K$ is nontrivial, we have that 
\[H_2 = \Gal(K(A[2^\infty])/F) \subsetneq \Gal(K(A[2^\infty])/K)  = \GSp_{2g}(\Z_2).\]
Further, $H(2) = \GSp_{2g}(\Z/2\Z)$ holds by the hypothesis of (\ref{A2:2}). Thus, applying Lemma \ref{ar-lem2} and Corollary \ref{sp-gsp-lift} to $A/F$, and observing that $H_2 = \Gal(F(A[2^\infty])/F)$, we find that $F/\Q$ is ramified at $2$. As $F$ is a subfield of $K(A_{\tors,(2)})$, this implies that $K(A_{\tors,(2)})/\Q$ is ramified at $2$. Thus 2 divides $\mathcal{B}_A$ by Corollary \ref{NOS-cor}.
\end{proof}

\section{Proof of Proposition \ref{prop-2} \label{prop-2-sec}}

We apply the constraints of \S\ref{constraints}  to prove Proposition \ref{prop-2}. Let $A$ and $g$ be as in the statement of the proposition. Let $\ell$ be an odd prime that divides $m_A$. By Lemmas \ref{ar-lem2}, \ref{gp-lem}, \ref{ar-lem1} and Corollary \ref{sp-gsp-lift}, we know
\[ \ell \text{ divides } \mathcal{B}_A \quad \text{or} \quad \ell \leq [\GSp_{2g}(\Z/\ell\Z): G(\ell)], \]
depending on whether $\Sp_{2g}(\Z/\ell\Z) \subseteq G(\ell)$ or not, respectively. Set $r' \coloneqq \rad'(m_A)$ and let $r'_{(2)}$ and ${\mathcal{B}_A}_{(2)}$ denote the odd-parts of $r'$ and $\mathcal{B}_A$, respectively (the odd part of an integer $n$ is $\frac{n}{2^k}$ where $2^k \mid\mid n$). Then,
\begin{align}
r'_{(2)} &\leq \prod_{\substack{\text{odd }\ell \mid m_A \\ \Sp_{2g}(\Z/\ell\Z) \subseteq G(\ell)}} \ell \prod_{\substack{\text{odd } \ell \mid m_A \nonumber\\ \Sp_{2g}(\Z/\ell\Z) \not\subseteq G(\ell)}} [\GSp_{2g}(\Z/\ell\Z) : G(\ell)]  \\ &\leq {\mathcal{B}_A}_{(2)} \cdot \left[\GSp_{2g}(\Z/r'_{(2)}\Z) : G(r'_{(2)}) \right]\!. \label{pf-prop42.1}
\end{align}
If $4 \nmid m_A$, then multiplying (\ref{pf-prop42.1}) through by 2, we obtain
\[ r' \leq 2 \cdot r'_{(2)} \leq 2 \cdot {\mathcal{B}_A}_{(2)} \cdot \left[\GSp_{2g}(\Z/r'_{(2)}\Z) : G(r'_{(2)}) \right] \leq 2 \cdot {\mathcal{B}_A} \cdot \left[\GSp_{2g}(\Z/r'\Z) : G(r') \right]\!. \]
If $4 \mid m_A$, then in particular $2 \mid m_A$, so by Lemma  \ref{ar-lem3}, we have that
\begin{align*}
2 \cdot [\GSp_{2g}(r'_{(2)}) : G(r'_{(2)})] \leq [\GSp_{2g}(r') : G(r')]   \quad \text{or} \quad   2 \cdot {\mathcal{B}_A}_{(2)} = \mathcal{B}_A.
\end{align*}
With this in mind, multiplying (\ref{pf-prop42.1}) through by $4$, we find that
\[r' = 4 r'_{(2)} \leq  4 \cdot{\mathcal{B}_A}_{(2)} \cdot  \left[\GSp_{2g}(\Z/r'_{(2)}\Z) : G(r'_{(2)}) \right]  \leq 2 \cdot \mathcal{B}_A \cdot \left[\GSp_{2g}(\Z/r'\Z) : G(r') \right]\!. \]
In either case, we see that the bound of Proposition \ref{prop-2} holds, completing its proof.

\section*{Acknowledgments}

The author thanks Nathan Jones for his valuable guidance. The author also thanks the anonymous referees for their comments that served to improve the paper.

\bibliographystyle{amsplain}
\bibliography{References}

\end{document}